\documentclass[12pt]{amsart}
\usepackage{geometry}   %????????
\usepackage[colorlinks,citecolor = red, linkcolor=blue,hyperindex]{hyperref}
\usepackage{euscript,verbatim, mathrsfs}
\usepackage[psamsfonts]{amssymb}
\usepackage{bbm}
\usepackage{graphicx}
\usepackage{float}
\usepackage{float, tikz}
\usepackage{extpfeil}
\usepackage[all, cmtip]{xy}
\usepackage{upref, xcolor, dsfont}
\usepackage{amsfonts,amsmath,amstext,amsbsy, amsopn,amsthm}
\usepackage{enumerate}
\usepackage{url}
\usepackage{mathtools}
\usepackage{bookmark}
\usepackage{euscript}
\usepackage{mathptmx}       % selects Times Roman as basic font
\usepackage{helvet}         % selects\textbf{\textbf{â¢}} Helvetica as sans-serif font
\usepackage{courier}        % selects Courier as typewriter font
\usepackage{type1cm}        % activate if the above 3 fonts are
\usepackage{multicol}       % used for the two-column index
\usepackage[bottom]{footmisc}% places footnotes at page bottom
\newtheorem{theorem}{Theorem}[section]
\newtheorem*{theorem*}{Theorem B}

\newtheorem{proposition}[theorem]{Proposition}

\newtheorem*{definition*}{Definition}
\newtheorem*{remark*}{Remark}

\newtheorem*{observation*}{Observation}

\newtheorem*{assumption*}{Assumption}
\newtheorem*{question*}{Question}
\newtheorem*{problem*}{Problem}

\newtheorem{mainthm}{Theorem}

\newtheorem{keylem}[mainthm]{Lemma}

\newcommand{\R}{\mathbb{R}}

\newcommand{\Var}{\mathrm{Var}}

\newcommand{\Conf}{\mathrm{Conf}}

\begin{document}
\title[Number rigid DPPs induced by generalized Cantor sets]{Number rigid determinantal point processes induced by generalized Cantor sets}

\author%[authorlabel1]
{Zhaofeng Lin}
\address%[authorlabel1]
{Zhaofeng Lin: School of Fundamental Physics and Mathematical Sciences, HIAS, University of Chinese Academy of Sciences, Hangzhou, 310024, China}
\email{linzhaofeng@ucas.ac.cn}

\author%[authorlabel1]
{Yanqi Qiu}
\address%[authorlabel1]
{Yanqi Qiu: School of Fundamental Physics and Mathematical Sciences, HIAS, University of Chinese Academy of Sciences, Hangzhou, 310024, China}
\email{yanqi.qiu@hotmail.com, yanqiqiu@ucas.ac.cn}

\author%[authorlabel1]
{Kai Wang}
\address%[authorlabel1]
{Kai Wang: School of Mathematical Sciences, Fudan University, Shanghai, 200433, China}
\email{kwang@fudan.edu.cn}

\begin{abstract}
We consider the Ghosh-Peres number rigidity of translation-invariant determinantal point processes on the real line $\mathbb{R}$, whose correlation kernels are induced by the Fourier transform of the indicators of generalized Cantor sets in the unit interval. Our main results show that for any given $\theta\in(0,1)$, there exists a generalized Cantor set with Lebesgue measure $\theta$, such that the corresponding determinantal point process is Ghosh-Peres number rigid.
\end{abstract}


\maketitle

\setcounter{equation}{0}

\section{Introduction and main results}
A point process is said to be Ghosh-Peres number rigid if for any compact subset in the phase space, the number of particles is almost surely determined by the restriction of the configuration to the complement of the compact subset (see below for its formal definition). In this paper, we consider the Ghosh-Peres number rigidity of translation-invariant determinantal point processes on the real line $\mathbb{R}$. Note that, in the case where the correlation kernel of the determinantal point process is Hermitian, a necessary condition for the number rigidity is that the correlation kernel is orthogonal projection. 

So we are led to consider the determinantal point processes  on $\R$ with correlation kernels:
\begin{align*}
	\widehat{\chi_M}(x-y)\quad x, y\in\mathbb{R},
\end{align*}
where $M$ are Borel subsets of $\mathbb{R}$ with finite Lebesgue measures. In the litterature, the number rigidity has  only be established in the case when the set $M$ is the union of a finite number of intervals. 

Therefore, a natural question is whether there exists number rigid determinantal point process on $\mathbb{R}$ induced by the union of an infinite number of intervals. The goal of the present paper is to answer affirmatively  this question.  Indeed, we construct some Ghosh-Peres number rigid determinantal point processes induced by the generalized Cantor sets, as well as their complementary sets which are union of an infinite number of intervals.

\subsection{Determinatal point processes}
We start by recalling some basic materials on point processes. Let $E$ be a locally compact Polish space. Denote by $\Conf(E)$ the space of all locally finite configurations over $E$, that is,
\begin{align*}
	\Conf(E)=\Big\{X=\sum_{i}\delta_{x_i}\,\Big|\,x_i\in E\,\,\text{and}\,\,X(\Delta)<\infty\,\,\text{for any compact subset}\,\,\Delta\subset E\Big\}.
\end{align*}
The configuration space $\Conf(E)$ is endowed with the vague topology, i.e., the weakest topology on $\Conf(E)$ such that for any compactly supported continuous function $f$ on $E$, the map $\Conf(E)\ni X\mapsto\int_{E}f\mathrm dX$ is continuous. It turns out that the configuration space $\Conf(E)$ equipped with the vague topology is a Polish space. By definition, a point process on $E$ is a Borel probability measure on $\Conf(E)$. For more details, see, e.g., \cite{DV1, DV2, Ka}.

A configuration $X=\sum_{i}\delta_{x_i}\in\Conf(E)$ is called simple if $X(\{x\})\leq1$ for all $x\in E$. In the simple case, we can identify $X=\sum_{i}\delta_{x_i}$ with the subset $X=\{x_i\}\subset E$. A point process $\mathbb{P}$ is called simple if $\mathbb{P}$-almost every $X\in\Conf(E)$ is simple.

Let $\mu$ be a Radon measure on $E$ and $K:E\times E\to\mathbb{C}$ be a measurable function. A simple point process $\mathbb{P}=\mathbb{P}_K$ on $E$ is called {\it determinantal} with correlation kernel $K$ with respect to the reference measure $\mu$, if for any $n\in\mathbb{N}_+$ and any 
bounded compactly supported measurable function $\varphi:E^n\to\mathbb{C}$,
\begin{align*}
	\mathbb{E}_{\mathbb{P}_K}\Big[\sum_{\substack{x_1,\cdots,\,x_k\in X\\\text{ distinct}}}\varphi(x_1,\cdots,x_n)\Big]=\int_{E^n}\varphi(x_1,\cdots,x_n)\det\big[K(x_i,x_j)\big]_{1\leq i,j\leq n}\mathrm{d}\mu(x_1)\cdots\mathrm{d}\mu(x_n).
\end{align*}
The reader is referred to \cite{Bo, HKPV, Le1, Le2, Le3, ST2, ST3} for more details on the theory of determinantal point processes.

Furthermore, we now consider the case that $K:E\times E\to\mathbb{C}$ is Hermitian. By slightly abusing the notation, we also denote $K$ the associated integral operator with the integral kernel $K(x,y)$, that is, $K:L^2(E,\mu)\to L^2(E,\mu)$ is defined by 
\begin{align*}
	Kf(x)=\int_E K(x,y)f(y)\mathrm{d}\mu(y),\quad x\in E,\quad f\in L^2(E,\mu).
\end{align*}
Suppose that the operator $K$ is locally trace class and positive contractive on $L^2(E,\mu)$ (here $K$ is locally trace class means that for each compact subset $\Delta\subset E$, the operator $\chi_{\Delta}K\chi _{\Delta}$ is trace class, where $\chi_{\Delta}$ denotes the multiplication operator by the indicator of $\Delta$), then by a theorem obtained by Macchi \cite{Ma} and Soshnikov \cite{So1}, as well as by Shirai and Takahashi \cite{ST1}, there exists a unique determinantal point process $\mathbb{P}_K$ on $E$ with correlation kernel $K$ with respect to the reference measure $\mu$.

For any bounded compactly supported measurable function $f:E\to\mathbb{C}$, define the {\it linear statistics}
\begin{align*}
	S_f(X)=\sum_{x\in X}f(x),\quad X\in\Conf(E).
\end{align*}
The moments of the linear statistics $S_f$ under the determinantal point process $\mathbb{P}_K$ can be calculated through the correlation kernel $K$. In particular, when $K$ is a reproducing kernel, the expectation and the variance of $S_f$ can be expressed respectively by
\begin{align*}
	\mathbb{E}_{\mathbb{P}_K}[S_f]&=\int_{E}f(x)K(x,x)\mathrm{d}\mu(x)
\end{align*}
and
\begin{align}\label{def-variance}
	\Var_{\mathbb{P}_K}(S_f)&=\frac{1}{2}\iint_{E^2}[f(x)-f(y)]^2|K(x,y)|^2\mathrm{d}\mu(x)\mathrm{d}\mu(y).
\end{align}
For more details, see, e.g., \cite{Gh, GP, So2}.

\subsection{Number rigid point processes}
We now recall the definition of the Ghosh-Peres number rigidity of point processes. Given a Borel subset $S\subset E$, let  ${\mathcal{F}}_S$ be the $\sigma$-algebra on $\Conf(E)$ generated by all functions $\Conf(E)\ni X\mapsto\#_B(X)=\#(B\cap X)$ with all Borel subsets $B\subset S$:
\begin{align*}
	{\mathcal{F}}_S=\sigma(\#_B: B\text{ is a Borel subset of } S).
\end{align*}
For any point process $\mathbb{P}$ on $E$, we denote by $\mathcal{F}_S^{\mathbb{P}}$ the completion of the $\sigma$-algebra ${\mathcal{F}}_S$ with respect to $\mathbb{P}$. A point process $\mathbb{P}$ on $E$ is called {\it number rigid} if for any compact subset $B\subset E$, the random variable $\#_B$ is ${\mathcal{F}}_{E\setminus B}^{\mathbb{P}}$-measurable. This definition of number rigidity is due to Ghosh \cite{Gh} where he showed that the sine-process is number rigid and Ghosh-Peres \cite{GP} where they showed that the Ginibre process and the zero set of Gaussian analytic function on the plane are number rigid.

Bufetov \cite{Bu} showed that the determinantal point processes with the Airy, the Bessel and the Gamma kernels are number rigid. Holroyd and Soo \cite{HS} showed that  the determinantal point process on the unit disc with the standard Bergman kernel is not number rigid. See also \cite{BQ} for an alternative proof of the non-rigid property of the Bergman point process. Qiu and Wang \cite{QW} constructed some non-trivial and natural number rigid determinantal point processes on the unit disc with sub-Bergman kernels in both deterministic and probabilistic methods. Bufetov, Dabrowski and Qiu \cite{BDQ} considered the linear rigidity of stochastic processes and obtained the linear rigidity for some determinantal point processes with translation-invariant kernels induced by a finite number of intervals of $\mathbb{R}$, which are of course number rigid. In addition, the authors gave a non-rigid determinantal point process on $\mathbb{R}^2$ in \cite{BDQ}. For more results on the number rigidity of point processes, we refer the reader to \cite{BFQ, BNQ, Gh1, GL, GL1, OS, RN}.

\subsection{Generalized Cantor sets}
Let $\alpha=\{\alpha_n\}_{n=1}^{\infty}$ be a series with $0<\alpha_n<1$, $n\in\mathbb{N}_+$. We will consider $\alpha$ as a sequence of ratios, and construct the generalized Cantor set by removing open intervals. The precise definition is given as follows.

Following the terminology in \cite{Fo}, if $I_0$ is a bounded interval and $\alpha_0\in(0,1)$, let us call the open interval with the same midpoint as $I_0$ and length equal to $\alpha_0$ times the length of $I_0$ the ``open middle $\alpha_0$-th'' of $I_0$.

For the sequence of ratios $\alpha=\{\alpha_n\}_{n=1}^{\infty}$, we can define a decreasing sequence $\{C_n\}_{n=1}^{\infty}$ of closed sets. Let us begin with the unit interval $C_0=[0,1]$, and for each $n\in\mathbb{N}_+$, $C_{n}$ is obtained by removing the open middle $\alpha_{n}$-th from each of the intervals that make up $C_{n-1}$. The resulting limiting set
\begin{align*}
	C=\bigcap_{n=1}^{\infty}C_{n}
\end{align*}
is called a {\it generalized Cantor set}. One can see that in the $n$-th iteration, the set we removed is composed of $2^{n-1}$ open intervals, denoted as $I_{n,k}$, $k=1,2,\cdots,2^{n-1}$. The  complementary set of the generalized Cantor set $C$ with respect to $[0,1]$ is
\begin{align*}
	I=[0,1]\setminus C=\bigsqcup_{n=1}^{\infty}\bigsqcup_{k=1}^{2^{n-1}}I_{n,k}.
\end{align*}

It is easy to show that the $2^{n-1}$ open intervals $I_{n,k}$, $k=1,2,\cdots,2^{n-1}$, have the same length
\begin{align}\label{def-length}
	l_n=\frac{(1-\alpha_1)(1-\alpha_2)\cdots(1-\alpha_{n-1})\alpha_n}{2^{n-1}}=\frac{1}{2^{n-1}}\Big[\prod_{k=1}^{n-1}(1-\alpha_k)-\prod_{k=1}^{n}(1-\alpha_k)\Big],
\end{align}
this is equivalent to
\begin{align}\label{def-ratio}
	\alpha_n=\frac{2^{n-1}l_n}{1-\sum_{k=1}^{n-1}2^{k-1}l_k}.
\end{align}
Hence the Lebesgue measure of $I$ is
\begin{align*}
	m(I)=\sum_{n=1}^{\infty}2^{n-1}l_n=1-\prod_{k=1}^{\infty}(1-\alpha_k).
\end{align*}
It follows that the Lebesgue measure of the generalized Cantor set $C$ is
\begin{align}\label{def-measureCantor}
	m(C)=\prod_{n=1}^{\infty}(1-\alpha_n)\in[0,1).
\end{align}

Note that the generalized Cantor set $C$ is compact, nowhere dense, totally disconnected, moreover, it has no isolated points. And the cardinality of $C$ is the cardinal number of the continuum. We refer the reader to \cite{Fo, Ro} for more details of the generalized Cantor sets.

\subsection{Main results}
For any sequence of ratios $\alpha=\{\alpha_n\}_{n=1}^{\infty}\subset(0,1)$, we have constructed the generalized Cantor set $C$ and its complementary set $I$ with respect to the unit interval $[0,1]$. Recall that the Fourier transform of a function $g\in L^1(\mathbb{R})$ is defined by
\begin{align*}
	\widehat{g}(\xi)=\int_{\mathbb{R}}g(x)e^{-2\pi ix\xi}\mathrm{d}x.
\end{align*}
Denote two translation-invariant kernels
\begin{align}\label{def-kernels}
	K_C(x,y)=\widehat{\chi_C}(x-y)\quad\text{and}\quad K_I(x,y)=\widehat{\chi_I}(x-y),\quad x,y\in\mathbb{R}.
\end{align}
By the theorem obtained by Macchi \cite{Ma} and Soshnikov \cite{So1}, as well as by Shirai and Takahashi \cite{ST1}, the kernels $K_C$ and $K_I$ induce two determinantal point processes on $\mathbb{R}$ with respect to the Lebesgue measure, denoted by $\mathbb{P}_{K_C}$ and $\mathbb{P}_{K_I}$ respectively.

Our first result is a sufficient condition with the length $l_n$ for the generalized Cantor set $C$ such that both the determinantal point processes $\mathbb{P}_{K_C}$ and $\mathbb{P}_{K_I}$ are number rigid.

\begin{theorem}\label{main-result1}
	For the generalized Cantor set $C$, if there exists $0<\delta<1$ such that
	\begin{align*}
		\sum_{n=1}^{\infty}\sqrt[4^{n/\delta}]{l_n}<\infty,
	\end{align*}
	then the determinantal point processes $\mathbb{P}_{K_C}$ and $\mathbb{P}_{K_I}$ are number rigid.
\end{theorem}

Our second result means that for any $\theta\in[0,1)$, one can find a generalized Cantor set $C$ with Lebesgue measure $\theta$, such that the corresponding determinantal point processes $\mathbb{P}_{K_C}$ and $\mathbb{P}_{K_I}$ are number rigid.

\begin{theorem}\label{main-result2}
	For any $\theta\in[0,1)$, there exists a generalized Cantor set $C$ satisfying $m(C)=\theta$ such that the determinantal point processes $\mathbb{P}_{K_C}$ and $\mathbb{P}_{K_I}$ are number rigid.
\end{theorem}

To prove our main results, we will use the classical method of Ghosh and Peres, i.e., the following Proposition~\ref{rigid-Ghosh-Peres}, which is a sufficient condition for number rigidity of a point process obtained by Ghosh \cite{Gh} and Ghosh-Peres \cite{GP}.

\begin{proposition}[Ghosh and Peres]\label{rigid-Ghosh-Peres}
	Let $\mathbb{P}$ be a point process on $E$. Assume that for any $\varepsilon>0$ and any compact subset $B\subset E$, there exists a bounded compactly supported measurable function $f:E\to\mathbb{C}$, such that $f\equiv1$ on $B$ and $\Var(S_f)<\varepsilon$. Then $\mathbb{P}$ is number rigid.
\end{proposition}

The proof of Theorem~\ref{main-result2} is constructive and is based on Theorem~\ref{main-result1}. The proof of Theorem~\ref{main-result1} relies on the following two key lemmas. 

\begin{keylem}\label{lem-Fourier}
	For the generalized Cantor set $C$, if there exists $0<\delta<1$ such that
	\begin{align*}
		\sum_{n=1}^{\infty}\sqrt[4^{n/\delta}]{l_n}<\infty,
	\end{align*}
	then there exists a constant $\lambda>0$ such that for any $\xi\neq0$,
	\begin{align*}
		|\widehat{\chi_C}(\xi)|^2\leq\lambda\frac{1+\log^{\delta}(1+|\xi|)}{|\xi|^2}\quad\text{and}\quad|\widehat{\chi_I}(\xi)|^2\leq\lambda\frac{1+\log^{\delta}(1+|\xi|)}{|\xi|^2}.
	\end{align*}
\end{keylem}

\begin{keylem}\label{lem-variance}
	For any $0<\varepsilon<1$, $0<\delta<1$ and $0<r<+\infty$, there exists a bounded compactly supported measurable function $\phi:\mathbb{R}\to\mathbb{R}$, such that $\phi\equiv1$ on $[-r,r]$ and
	\begin{align*}
		\iint_{\mathbb{R}^2}\Big[\frac{\phi(x)-\phi(y)}{x-y}\Big]^2\big[1+\log^{\delta}(1+|x-y|)\big]\mathrm{d}x\mathrm{d}y<\varepsilon.
	\end{align*}
\end{keylem}

We will prove Lemmas~\ref{lem-Fourier} and \ref{lem-variance} respectively in Sections~\ref{proof-lem-Fourier} and \ref{proof-lem-variance}. The proof of Lemma~\ref{lem-variance} is constructive inspited by Bufetov \cite{Bu}.

We end this section by providing the standard derivation of Theorem~\ref{main-result1} from Lemmas~\ref{lem-Fourier} and \ref{lem-variance}, and the constructive derivation of  Theorem~\ref{main-result2} from Theorem~\ref{main-result1}.

{\flushleft \it The derivation of Theorem~\ref{main-result1} from Lemmas~\ref{lem-Fourier} and \ref{lem-variance}.}
By the definitions of $K_C$ and $K_I$ in \eqref{def-kernels}, it follows from Lemma~\ref{lem-Fourier} that for any $x,y\in\mathbb{R}$, $x\neq y$, we have
\begin{align*}
	|K_C(x,y)|^2\leq\lambda\frac{1+\log^{\delta}(1+|x-y|)}{|x-y|^2}\quad\text{and}\quad|K_I(x,y)|^2\leq\lambda\frac{1+\log^{\delta}(1+|x-y|)}{|x-y|^2}.
\end{align*}
Hence the number rigidity of the determinantal point processes $\mathbb{P}_{K_C}$ and $\mathbb{P}_{K_I}$ can be obtained immediately by formula \eqref{def-variance}, Proposition~\ref{rigid-Ghosh-Peres} and Lemma~\ref{lem-variance}.

{\flushleft \it The derivation of Theorem~\ref{main-result2} from Theorem~\ref{main-result1}.}
Let $\{u_n\}_{n=1}^{\infty}$ be any sequence of positive number satisfying $\sum_{n=1}^{\infty}u_n<\infty$. Since $\lim\limits_{n\to\infty}u_n=0$, for any $0<\delta<1$, there exists a constant $c>0$ such that $2^{n-1}u_n^{4^{n/\delta}}\leq(2u_n)^{4^{n/\delta}}\leq cu_n$ for any $n\in\mathbb{N}_+$, and hence $\sum_{n=1}^{\infty}2^{n-1}u_n^{4^{n/\delta}}<\infty$. Now for any $\theta\in[0,1)$, let $\Theta>0$ be the constant such that
\begin{align*}
	\Theta\sum_{n=1}^{\infty}2^{n-1}u_n^{4^{n/\delta}}=1-\theta.
\end{align*}
By \eqref{def-length} and \eqref{def-ratio}, for each $n\in\mathbb{N}_+$, denote
\begin{align*}
	l_n=\Theta u_n^{4^{n/\delta}}\quad\text{and}\quad\alpha_n=\frac{2^{n-1}l_n}{1-\sum_{k=1}^{n-1}2^{k-1}l_k}.
\end{align*}
Consider the generalized Cantor set $C$ corresponding to $\{\alpha_n\}_{n=1}^{\infty}$, we know that the length $l_n$ satisfies
\begin{align*}
	\sum_{n=1}^{\infty}\sqrt[4^{n/\delta}]{l_n}=\sum_{n=1}^{\infty}\sqrt[4^{n/\delta}]{\Theta}\,u_n\leq(1+\Theta)\sum_{n=1}^{\infty}u_n<\infty.
\end{align*}
Thus by Theorem~\ref{main-result1}, the corresponding  determinantal point processes $\mathbb{P}_{K_C}$ and $\mathbb{P}_{K_I}$ are number rigid. And by \eqref{def-measureCantor}, we have that
\begin{align*}
	m(C)=\prod_{n=1}^{\infty}(1-\alpha_n)=\prod_{n=1}^{\infty}\frac{1-\sum_{k=1}^{n}2^{k-1}l_k}{1-\sum_{k=1}^{n-1}2^{k-1}l_k}=1-\sum_{k=1}^{\infty}2^{k-1}l_k=1-\Theta\sum_{k=1}^{\infty}2^{k-1}u_k^{4^{k/\delta}}=\theta.
\end{align*}
This clearly completes the proof of Theorem~\ref{main-result2}.

{\flushleft\bf Acknowledgements.}
YQ's work is supported by the National Natural Science Foundation of China (No.12288201).
KW's work is supported by the National Natural Science Foundation of China (No.12231005, No.12326376) and the Shanghai Technology Innovation Project (21JC1400800).

\section{The proof of Lemma~\ref{lem-Fourier}}\label{proof-lem-Fourier}
In this section, we are going to prove Lemma~\ref{lem-Fourier}. We first estimate $|\widehat{\chi_I}|^2$ under the condition $\sum_{n=1}^{\infty}\sqrt[4^{n/\delta}]{l_n}<\infty$ with $0<\delta<1$. The estimate of $|\widehat{\chi_C}|^2$ follows from that of $|\widehat{\chi_I}|^2$.

\subsection{The estimate of $|\widehat{\chi_I}|^2$}
Recall that the complementary set of the generalized Cantor set $C$ with respect to $[0,1]$ is given by
\begin{align*}
	I=\bigsqcup_{n=1}^{\infty}\bigsqcup_{k=1}^{2^{n-1}}I_{n,k}.
\end{align*}
For each $n\in\mathbb{N}_+$ and $k=1,2,\cdots,2^{n-1}$, we have
\begin{align*}
	\widehat{\chi}_{I_{n,k}}(\xi)=\int_{\mathbb{R}}\chi_{I_{n,k}}(x)e^{-2\pi ix\xi}\mathrm{d}x=\frac{e^{-2\pi ir(I_{n,k})\xi}\big(1-e^{-2\pi il_n\xi}\big)}{2\pi i\xi},\quad\xi\neq0,
\end{align*}
where $r(I_{n,k})$ is the right endpoint of the interval $I_{n,k}$. It follows that
\begin{align*}
	|\widehat{\chi}_{I_{n,k}}(\xi)|=\frac{|\sin(\pi l_n\xi)|}{\pi|\xi|},\quad\xi\neq0.
\end{align*}
Thus for any $\xi\neq0$, we have
\begin{align*}
	|\widehat{\chi}_{I}(\xi)|=\Big|\sum_{n=1}^{\infty}\sum_{k=1}^{2^{n-1}}\widehat{\chi}_{I_{n,k}}(\xi)\Big|\leq\sum_{n=1}^{\infty}\sum_{k=1}^{2^{n-1}}|\widehat{\chi}_{I_{n,k}}(\xi)|=\frac{1}{\pi|\xi|}\sum_{n=1}^{\infty}2^{n-1}|\sin(\pi l_n\xi)|.
\end{align*}
We will show that there exists a constant $\kappa>0$ such that 
\begin{align*}
	\sum_{n=1}^{\infty}2^{n-1}|\sin(\pi l_n\xi)|\leq\kappa\big[1+\log^{\frac{\delta}{2}}(1+|\xi|)\big],
\end{align*}
so that there exists a constant $\lambda>0$ such that for any $\xi\neq0$,
\begin{align*}
	|\widehat{\chi_I}(\xi)|^2\leq\lambda\frac{1+\log^{\delta}(1+|\xi|)}{|\xi|^2}.
\end{align*}

Indeed, set
\begin{align*}
	a(\xi)=\frac{\log\big(1+A\big[1+\log^{\frac{\delta}{2}}(1+|\xi|)\big]\big)}{\log2},
\end{align*}
here $A>0$ is a fixed sufficiently large constant. Let us write the following decomposition
\begin{align*}
	\sum_{n=1}^{\infty}2^{n-1}|\sin(\pi l_n\xi)|=\sum_{n=1}^{\lfloor a(\xi)\rfloor}2^{n-1}|\sin(\pi l_n\xi)|+\sum_{n=1+\lfloor a(\xi)\rfloor}^{\infty}2^{n-1}|\sin(\pi l_n\xi)|,
\end{align*}
where $\lfloor a(\xi)\rfloor$ is the biggest integer that does not exceed $a(\xi)$. It is easy to see that $\frac{n}{4^{n/\delta}}\leq1$ for any $n\in\mathbb{N}_+$. Then since $|\sin(\pi l_n\xi)|\leq1$,
\begin{align*}
	|\sin(\pi l_n\xi)|\leq|\sin(\pi l_n\xi)|^{\frac{n}{4^{n/\delta}}},
\end{align*}
hence
\begin{align*}
	\sum_{n=1}^{\infty}2^{n-1}|\sin(\pi l_n\xi)|\leq\sum_{n=1}^{\lfloor a(\xi)\rfloor}2^{n-1}+\sum_{n=1+\lfloor a(\xi)\rfloor}^{\infty}2^{n-1}|\sin(\pi l_n\xi)|^{\frac{n}{4^{n/\delta}}}.
\end{align*}
Noticing that
\begin{align*}
	\sum_{n=1}^{\lfloor a(\xi)\rfloor}2^{n-1}=2^{\lfloor a(\xi)\rfloor}-1\leq2^{a(\xi)}-1=A\big[1+\log^{\frac{\delta}{2}}(1+|\xi|)\big]
\end{align*}
and
\begin{align*}
	|\sin(\pi l_n\xi)|\leq\pi l_n|\xi|\leq\pi l_n(1+|\xi|),
\end{align*}
we obtain
\begin{align*}
	&\quad\,\,\sum_{n=1}^{\infty}2^{n-1}|\sin(\pi l_n\xi)|\\
	&\leq A\big[1+\log^{\frac{\delta}{2}}(1+|\xi|)\big]+\sum_{n=1+\lfloor a(\xi)\rfloor}^{\infty}2^{n-1}\pi^{\frac{n}{4^{n/\delta}}}l_n^{\frac{n}{4^{n/\delta}}}(1+|\xi|)^{\frac{n}{4^{n/\delta}}}\\
	&\leq A\big[1+\log^{\frac{\delta}{2}}(1+|\xi|)\big]+\frac{\pi}{2}\sum_{n=1+\lfloor a(\xi)\rfloor}^{\infty}\big(2\sqrt[4^{n/\delta}]{l_n}\big)^n(1+|\xi|)^{\frac{n}{4^{n/\delta}}}.
\end{align*}
Since the constant $A>0$ is sufficiently large, we know that for any $n>a(\xi)$,
\begin{align*}
	(1+|\xi|)^{\frac{n}{4^{n/\delta}}}\leq(1+|\xi|)^{\frac{a(\xi)}{4^{a(\xi)/\delta}}}.
\end{align*}
By a direct calculation,
\begin{align*}
	(1+|\xi|)^{\frac{a(\xi)}{4^{a(\xi)/\delta}}}=\exp\left(\frac{\log(1+|\xi|)\log\big(1+A\big[1+\log^{\frac{\delta}{2}}(1+|\xi|)\big]\big)}{\big(1+A\big[1+\log^{\frac{\delta}{2}}(1+|\xi|)\big]\big)^{\frac{2}{\delta}}\log2}\right),
\end{align*}
which implies that
\begin{align*}
	(1+|\xi|)^{\frac{a(\xi)}{4^{a(\xi)/\delta}}}\leq\big(1+A\big[1+\log^{\frac{\delta}{2}}(1+|\xi|)\big]\big)^{\frac{1}{A^{\frac{\delta}{2}}\log2}}\leq(A+1)\big[1+\log^{\frac{\delta}{2}}(1+|\xi|)\big].
\end{align*}
Hence for any $n>a(\xi)$,
\begin{align*}
	(1+|\xi|)^{\frac{n}{4^{n/\delta}}}\leq(A+1)\big[1+\log^{\frac{\delta}{2}}(1+|\xi|)\big].
\end{align*}
Therefore,
\begin{align*}
	&\quad\,\,\sum_{n=1}^{\infty}2^{n-1}|\sin(\pi l_n\xi)|\\
	&\leq A\big[1+\log^{\frac{\delta}{2}}(1+|\xi|)\big]+\frac{\pi}{2}\sum_{n=1+\lfloor a(\xi)\rfloor}^{\infty}\big(2\sqrt[4^{n/\delta}]{l_n}\big)^n(A+1)\big[1+\log^{\frac{\delta}{2}}(1+|\xi|)\big]\\
	&\leq\Big[A+\frac{\pi(A+1)}{2}\sum_{n=1}^{\infty}\big(2\sqrt[4^{n/\delta}]{l_n}\big)^n\Big]\big[1+\log^{\frac{\delta}{2}}(1+|\xi|)\big].
\end{align*}

Since
\begin{align*}
	\sum_{n=1}^{\infty}\sqrt[4^{n/\delta}]{l_n}<\infty,
\end{align*}
we have $\lim\limits_{n\to\infty}\sqrt[4^{n/\delta}]{l_n}=0$. It follows that $\big(2\sqrt[4^{n/\delta}]{l_n}\big)^n\leq2\sqrt[4^{n/\delta}]{l_n}$ for sufficiently large $n$, hence
\begin{align*}
	\sum_{n=1}^{\infty}\big(2\sqrt[4^{n/\delta}]{l_n}\big)^n\leq2\sum_{n=1}^{\infty}\sqrt[4^{n/\delta}]{l_n}<\infty.
\end{align*}
Denote
\begin{align*}
	\kappa=A+\frac{\pi(A+1)}{2}\sum_{n=1}^{\infty}\big(2\sqrt[4^{n/\delta}]{l_n}\big)^n,
\end{align*}
then
\begin{align*}
	\sum_{n=1}^{\infty}2^{n-1}|\sin(\pi l_n\xi)|\leq\kappa\big[1+\log^{\frac{\delta}{2}}(1+|\xi|)\big].
\end{align*}
Consequently, there exists a constant $\lambda>0$ such that for any $\xi\neq0$,
\begin{align}\label{estimate-widehatchiI}
	|\widehat{\chi_I}(\xi)|^2\leq\lambda\frac{1+\log^{\delta}(1+|\xi|)}{|\xi|^2}.
\end{align}

\subsection{The estimate of $|\widehat{\chi_C}|^2$}
Rewrite the generalized Cantor set in the form of 
\begin{align*}
	C=[0,1]\setminus I,
\end{align*}
then
\begin{align*}
	\widehat{\chi_C}=\widehat{\chi_{[0,1]}}-\widehat{\chi_I}
\end{align*}
and hence
\begin{align*}
	|\widehat{\chi_C}|^2\leq|\widehat{\chi_{[0,1]}}|^2+|\widehat{\chi_I}|^2+2|\widehat{\chi_{[0,1]}}||\widehat{\chi_I}|.
\end{align*}
Note that for any $\xi\neq0$,
\begin{align*}
	|\widehat{\chi_{[0,1]}}(\xi)|=\Big|\frac{\sin(\pi\xi)}{\pi\xi}\Big|\leq\frac{1}{\pi|\xi|}.
\end{align*}
Combine with the estimate of $|\widehat{\chi_I}|$, i.e., the estimate of $|\widehat{\chi_I}|^2$ obtained in \eqref{estimate-widehatchiI}, we conclude that there exists a constant $\lambda>0$ such that for any $\xi\neq0$,
\begin{align}\label{estimate-widehatchiC}
	|\widehat{\chi_C}(\xi)|^2\leq\lambda\frac{1+\log^{\delta}(1+|\xi|)}{|\xi|^2}.
\end{align}

\subsection{The conclusion of Lemma~\ref{lem-Fourier}}
Lemma~\ref{lem-Fourier} follows from \eqref{estimate-widehatchiI} and \eqref{estimate-widehatchiC}.

\section{The proof of Lemma~\ref{lem-variance}}\label{proof-lem-variance}
In this section, we are going to prove Lemma~\ref{lem-variance}. For any $0<r<R<+\infty$, denote
\begin{align*}
	\phi^{(r,R)}(x)=\left\{\begin{array}{ll}
		1,\qquad\qquad\qquad\,\,\,\,|x|\leq r,\\
		1-\frac{\log(1+|x|-r)}{\log(1+R-r)},\quad r<|x|<R,\\
		0,\qquad\qquad\qquad\,\,\,\,|x|\geq R.
	\end{array}\right.
\end{align*}
The function $\phi^{(r,R)}$ is inspited by Bufetov \cite{Bu}, which is a bounded compactly supported measurable function such that $\phi^{(r,R)}\equiv1$ on $[-r,r]$. We will show that for any $0<\delta<1$ and $0<r<+\infty$,
\begin{align}\label{estimate-variance}
	\lim_{R\to+\infty}\iint_{\mathbb{R}^2}\Big[\frac{\phi^{(r,R)}(x)-\phi^{(r,R)}(y)}{x-y}\Big]^2\big[1+\log^{\delta}(1+|x-y|)\big]\mathrm{d}x\mathrm{d}y=0.
\end{align}
When $|x|,|y|\leq r$ or $|x|,|y|\geq R$, the expression of the integral is equal to zero. So we shall estimate our integral over the following four domains:
\begin{align*}
	D_1(r,R)=\{(x,y)\in\mathbb{R}^2:r<|x|,|y|<R\};\quad D_2(r,R)=\{(x,y)\in\mathbb{R}^2:|x|<r<|y|<R\};
\end{align*}
\begin{align*}
	D_3(r,R)=\{(x,y)\in\mathbb{R}^2:|x|<r<R<|y|\};\quad D_4(r,R)=\{(x,y)\in\mathbb{R}^2:r<|x|<R<|y|\}.
\end{align*}

\subsection{The estimate over $D_1(r,R)$}
The integral over the domain $D_1(r,R)$ can be expressed explicitly as
\begin{align*}
	J_{R,1}&=\frac{1}{\log^2(1+R-r)}\int_{r<|x|<R}\int_{r<|y|<R}\Big[\frac{\log(1+|x|-r)-\log(1+|y|-r)}{x-y}\Big]^2\\
	&\qquad\qquad\qquad\qquad\qquad\qquad\qquad\times\big[1+\log^{\delta}(1+|x-y|)\big]\mathrm{d}x\mathrm{d}y.
\end{align*}
Noticing that there exists a constant $c(r)>0$ such that for any $|x|,|y|>r$,
\begin{align*}
	\big|\log(1+|x|-r)-\log(1+|y|-r)\big|\leq c(r)\big|\log|x|-\log|y|\big|,
\end{align*}
and by using
\begin{align*}
	|x-y|\geq\big||x|-|y|\big|,
\end{align*}
we have
\begin{align*}
	J_{R,1}&\leq\frac{c^2(r)\big[1+\log^{\delta}(1+2R)\big]}{\log^2(1+R-r)}\int_{r<|x|<R}\int_{r<|y|<R}\Big(\frac{\log|x|-\log|y|}{|x|-|y|}\Big)^2\mathrm{d}x\mathrm{d}y\\
	&=\frac{4c^2(r)\big[1+\log^{\delta}(1+2R)\big]}{\log^2(1+R-r)}\int_{r}^{R}\int_{r}^{R}\Big(\frac{\log x-\log y}{x-y}\Big)^2\mathrm{d}x\mathrm{d}y\\
	&=\frac{4c^2(r)\big[1+\log^{\delta}(1+2R)\big]}{\log^2(1+R-r)}\int_{r}^{R}\int_{r}^{R}\Big(\frac{\log\frac{y}{x}}{\frac{y}{x}-1}\Big)^2\frac{1}{x^2}\mathrm{d}x\mathrm{d}y.
\end{align*}
By a substitution $t=\frac{y}{x}$, we obtain
\begin{align*}
	J_{R,1}&\leq\frac{4c^2(r)\big[1+\log^{\delta}(1+2R)\big]}{\log^2(1+R-r)}\int_{r}^{R}\int_{\frac{r}{x}}^{\frac{R}{x}}\Big(\frac{\log t}{t-1}\Big)^2\frac{1}{x}\mathrm{d}t\mathrm{d}x\\
	&\leq\frac{4c^2(r)\big[1+\log^{\delta}(1+2R)\big]}{\log^2(1+R-r)}\int_{r}^{R}\frac{1}{x}\mathrm{d}x\int_{\frac{r}{R}}^{\frac{R}{r}}\Big(\frac{\log t}{t-1}\Big)^2\mathrm{d}t\\
	&=\frac{4c^2(r)\big[1+\log^{\delta}(1+2R)\big](\log R-\log r)}{\log^2(1+R-r)}\int_{\frac{r}{R}}^{\frac{R}{r}}\Big(\frac{\log t}{t-1}\Big)^2\mathrm{d}t.
\end{align*}
Since
\begin{align*}
	\int_{0}^{+\infty}\Big(\frac{\log t}{t-1}\Big)^2\mathrm{d}t<\infty,
\end{align*}
we get
\begin{align}\label{estimate-JR1}
	J_{R,1}\leq\frac{4c^2(r)\big[1+\log^{\delta}(1+2R)\big](\log R-\log r)}{\log^2(1+R-r)}\int_{0}^{+\infty}\Big(\frac{\log t}{t-1}\Big)^2\mathrm{d}t=O\Big(\frac{1}{\log^{1-\delta}R}\Big).
\end{align}

\subsection{The estimate over $D_2(r,R)$}
The integral over the domain $D_2(r,R)$ can be expressed explicitly as
\begin{align*}
	J_{R,2}=\frac{1}{\log^2(1+R-r)}\int_{|x|<r}\int_{r<|y|<R}\Big[\frac{\log(1+|y|-r)}{x-y}\Big]^2\big[1+\log^{\delta}(1+|x-y|)\big]\mathrm{d}x\mathrm{d}y.
\end{align*}
Noticing that when $|x|<r<|y|$, we have
\begin{align*}
	|x-y|\geq|y|-r,
\end{align*}
thus
\begin{align*}
	J_{R,2}&\leq\frac{2r\big[1+\log^{\delta}(1+R+r)\big]}{\log^2(1+R-r)}\int_{r<|y|<R}\Big[\frac{\log(1+|y|-r)}{|y|-r}\Big]^2\mathrm{d}y\\
	&=\frac{4r\big[1+\log^{\delta}(1+R+r)\big]}{\log^2(1+R-r)}\int_{r}^{R}\Big[\frac{\log(1+y-r)}{y-r}\Big]^2\mathrm{d}y\\
	&=\frac{4r\big[1+\log^{\delta}(1+R+r)\big]}{\log^2(1+R-r)}\int_{0}^{R-r}\Big[\frac{\log(1+t)}{t}\Big]^2\mathrm{d}t.
\end{align*}
Since
\begin{align*}
	\int_{0}^{+\infty}\Big[\frac{\log(1+t)}{t}\Big]^2\mathrm{d}t<\infty,
\end{align*}
we get
\begin{align}\label{estimate-JR2}
	J_{R,2}\leq\frac{4r\big[1+\log^{\delta}(1+R+r)\big]}{\log^2(1+R-r)}\int_{0}^{+\infty}\Big[\frac{\log(1+t)}{t}\Big]^2\mathrm{d}t=O\Big(\frac{1}{\log^{2-\delta}R}\Big).
\end{align}

\subsection{The estimate over $D_3(r,R)$}
The integral over the domain $D_3(r,R)$ can be expressed explicitly as
\begin{align*}
	J_{R,3}=\int_{|x|<r}\int_{|y|>R}\frac{1}{(x-y)^2}\big[1+\log^{\delta}(1+|x-y|)\big]\mathrm{d}x\mathrm{d}y.
\end{align*}
Noticing that
\begin{align*}
	\log(1+|x-y|)\leq|x-y|,
\end{align*}
and when $|x|<r<R<|y|$, we have
\begin{align*}
	|x-y|\geq|y|-r,
\end{align*}
thus
\begin{align*}
	J_{R,3}&\leq\int_{|x|<r}\int_{|y|>R}\Big[\frac{1}{(x-y)^2}+\frac{1}{|x-y|^{2-\delta}}\Big]\mathrm{d}x\mathrm{d}y\\
	&\leq2r\int_{|y|>R}\Big[\frac{1}{(|y|-r)^2}+\frac{1}{(|y|-r)^{2-\delta}}\Big]\mathrm{d}y\\
	&=4r\int_{R}^{+\infty}\Big[\frac{1}{(y-r)^2}+\frac{1}{(y-r)^{2-\delta}}\Big]\mathrm{d}y.
\end{align*}
Since
\begin{align*}
	\int_{R}^{+\infty}\Big[\frac{1}{(y-r)^2}+\frac{1}{(y-r)^{2-\delta}}\Big]\mathrm{d}y=\frac{1}{R-r}+\frac{1}{1-\delta}\frac{1}{(R-r)^{1-\delta}},
\end{align*}
we get
\begin{align}\label{estimate-JR3}
	J_{R,3}\leq4r\Big[\frac{1}{R-r}+\frac{1}{1-\delta}\frac{1}{(R-r)^{1-\delta}}\Big]=O\Big(\frac{1}{R^{1-\delta}}\Big).
\end{align}

\subsection{The estimate over $D_4(r,R)$}
The integral over the domain $D_4(r,R)$ can be expressed explicitly as
\begin{align*}
	J_{R,4}&=\frac{1}{\log^2(1+R-r)}\int_{r<|x|<R}\int_{|y|>R}\Big[\frac{\log(1+R-r)-\log(1+|x|-r)}{x-y}\Big]^2\\
	&\qquad\qquad\qquad\qquad\qquad\qquad\quad\times\big[1+\log^{\delta}(1+|x-y|)\big]\mathrm{d}x\mathrm{d}y.
\end{align*}
Noticing that in the case $r<|x|<R<|y|$, we have
\begin{align*}
	&\quad\,\,\big|\log(1+R-r)-\log(1+|x|-r)\big|=\log(1+R-r)-\log(1+|x|-r)\\
	&\leq\log(1+|y|-r)-\log(1+|x|-r)=\big|\log(1+|y|-r)-\log(1+|x|-r)\big|\\
	&\leq c(r)\big|\log|x|-\log|y|\big|,
\end{align*}
and by using
\begin{align*}
	|x-y|\geq\big||x|-|y|\big|,
\end{align*}
we have
\begin{align*}
	J_{R,4}&\leq\frac{c^2(r)}{\log^2(1+R-r)}\int_{r<|x|<R}\int_{|y|>R}\Big(\frac{\log|x|-\log|y|}{|x|-|y|}\Big)^2\big[1+\log^{\delta}(1+|x|+|y|)\big]\mathrm{d}x\mathrm{d}y\\
	&=\frac{4c^2(r)}{\log^2(1+R-r)}\int_{r}^{R}\int_{R}^{+\infty}\Big(\frac{\log x-\log y}{x-y}\Big)^2\big[1+\log^{\delta}(1+x+y)\big]\mathrm{d}x\mathrm{d}y\\
	&=\frac{4c^2(r)}{\log^2(1+R-r)}\int_{r}^{R}\int_{R}^{+\infty}\Big(\frac{\log\frac{y}{x}}{\frac{y}{x}-1}\Big)^2\big[1+\log^{\delta}\big(1+x(1+\frac{y}{x})\big)\big]\frac{1}{x^2}\mathrm{d}x\mathrm{d}y.
\end{align*}
By a substitution $t=\frac{y}{x}$, we obtain
\begin{align*}
	J_{R,4}&\leq\frac{4c^2(r)}{\log^2(1+R-r)}\int_{r}^{R}\int_{\frac{R}{x}}^{+\infty}\Big(\frac{\log t}{t-1}\Big)^2\big[1+\log^{\delta}\big(1+x(1+t)\big)\big]\frac{1}{x}\mathrm{d}t\mathrm{d}x\\
	&\leq\frac{4c^2(r)}{\log^2(1+R-r)}\int_{r}^{R}\frac{1}{x}\mathrm{d}x\int_{1}^{+\infty}\Big(\frac{\log t}{t-1}\Big)^2\big[1+\log^{\delta}\big(1+R(1+t)\big)\big]\mathrm{d}t\\
	&=\frac{4c^2(r)(\log R-\log r)}{\log^2(1+R-r)}\int_{1}^{+\infty}\Big(\frac{\log t}{t-1}\Big)^2\big[1+\log^{\delta}\big(1+R(1+t)\big)\big]\mathrm{d}t.
\end{align*}
Note that
\begin{align*}
	&\quad\,\,\log^{\delta}\big(1+R(1+t)\big)\leq\log^{\delta}\big((1+R)(1+t)\big)=\big[\log(1+R)+\log(1+t)\big]^{\delta}\\
	&=\frac{\log(1+R)+\log(1+t)}{\big[\log(1+R)+\log(1+t)\big]^{1-\delta}}\leq\frac{\log(1+R)+\log(1+t)}{\log^{1-\delta}(1+R)}=\log^{\delta}(1+R)+\frac{\log(1+t)}{\log^{1-\delta}(1+R)},
\end{align*}
hence
\begin{align*}
	J_{R,4}\leq\frac{4c^2(r)(\log R-\log r)}{\log^2(1+R-r)}\int_{1}^{+\infty}\Big(\frac{\log t}{t-1}\Big)^2\Big[1+\log^{\delta}(1+R)+\frac{\log(1+t)}{\log^{1-\delta}(1+R)}\Big]\mathrm{d}t.
\end{align*}
Since
\begin{align*}
	\int_{1}^{+\infty}\Big(\frac{\log t}{t-1}\Big)^2\mathrm{d}t<\infty\quad\text{and}\quad\int_{1}^{+\infty}\Big(\frac{\log t}{t-1}\Big)^2\log(1+t)\mathrm{d}t<\infty,
\end{align*}
we get
\begin{align}\label{estimate-JR4}
	J_{R,4}\leq O\Big(\frac{1}{\log R}\Big)+O\Big(\frac{1}{\log^{1-\delta}R}\Big)+O\Big(\frac{1}{\log^{2-\delta}R}\Big)=O\Big(\frac{1}{\log^{1-\delta}R}\Big).
\end{align}

\subsection{The conclusion of Lemma~\ref{lem-variance}}
Based on the estimates \eqref{estimate-JR1}-\eqref{estimate-JR4}, we obtain the limit in \eqref{estimate-variance}. Then Lemma~\ref{lem-variance} follows from \eqref{estimate-variance}.

\end{document}